\documentclass[11pt]{amsart}

\usepackage{amsmath}
\usepackage{amssymb}

\newtheorem{theorem}{Theorem}[section]
\newtheorem{claim}[theorem]{Claim}

\newtheorem{corollary}[theorem]{Corollary}

\theoremstyle{definition}
\newtheorem{definition}[theorem]{Definition}

\newtheorem{question}[theorem]{Question}

\theoremstyle{remark}
\newtheorem{remark}[theorem]{Remark}

\newcount\skewfactor
\def\mathunderaccent#1#2 {\let\theaccent#1\skewfactor#2
\mathpalette\putaccentunder}
\def\putaccentunder#1#2{\oalign{$#1#2$\crcr\hidewidth
\vbox to.2ex{\hbox{$#1\skew\skewfactor\theaccent{}$}\vss}\hidewidth}}

\def\smallbox#1{\leavevmode\thinspace\hbox{\vrule\vtop{\vbox
   {\hrule\kern1pt\hbox{\vphantom{\tt/}\thinspace{\tt#1}\thinspace}}
   \kern1pt\hrule}\vrule}\thinspace}

\newcommand{\cf}{{\rm cf}}

\def\qedref#1{$\qed_{\reforiginal{#1}}$}

\title{Dense free sets}
\author{Shimon Garti}
\address{Institute of Mathematics,
 The Hebrew University of Jerusalem,
 Jerusalem 91904, Israel}
\email{shimon.garty@mail.huji.ac.il}
\thanks{}

\subjclass[2010]{Primary: 54A25, Secondary: 03E02}
\keywords{Free sets, reaping number, everywhere dense set, polarized relation}

\begin{document}
\let\labeloriginal\label
\let\reforiginal\ref

\begin{abstract}
Let $\kappa=2^\omega$, and assume $f:\mathbb{R}\rightarrow\mathcal{P}(\mathbb{R})$ satisfies the intersection properties $C(\omega,\kappa)$ and $C(\kappa,\omega)$. We prove that if $\mathfrak{r}<\cf(\kappa)$ then there exists a dense free set for $f$.
\end{abstract}

\maketitle

\newpage

\section{Introduction}

The basic notion of this paper is the following:

\begin{definition}
\label{ffrree} Free sets. \newline
Let $f$ be a function from $\theta$ into $\mathcal{P}(\theta)$. \newline
A set $A\subseteq\theta$ is free (for $f$) iff $x\notin f(y)$ whenever $\{x,y\}\subseteq A$.
\end{definition}

Free sets are quite useful in many branches of combinatorial set theory, and the basic problems are the existence of large free sets. Two simple examples of functions show that one must make some assumptions on $f$ in order to get an infinite free set. The following definition and the examples below are phrased in \cite{MR3310341}:

\begin{definition}
\label{iintersect} The intersection property $C(\lambda,\mu)$. \newline
Let $f$ be a function from $S$ into $\mathcal{P}(S)$. \newline
We say that $f$ satisfies the property $C(\lambda,\mu)$ iff $|\bigcap\{f(x):x\in T\}|<\mu$ for every subset $T$ of $S$ of size $\lambda$.
\end{definition}

Let $\kappa$ be an infinite cardinal. The initial-segment coloring $f(\alpha)=\{\beta:\beta<\alpha\}$ satisfies $C(\omega,\kappa)$ but not $C(\kappa,\omega)$. Clearly, it has no infinite free set (actually, not even two-element free set). The end-segment coloring $f(\alpha)=\{\beta:\beta>\alpha\}$ satisfies $C(\kappa,\omega)$ but not $C(\omega,\kappa)$. This function fails similarly to have an infinite free set. Consequently, one has to assume the properties $C(\kappa,\omega)$ and $C(\omega,\kappa)$ in order to exclude trivial cases:

\begin{definition}
\label{rreasonable} Reasonable set mappings. \newline
A function $f:\kappa\rightarrow\mathcal{P}(\kappa)$ is called $\kappa$-reasonable iff $f$ satisfies both $C(\kappa,\omega)$ and $C(\omega,\kappa)$.
\end{definition}

However, the above trivial restriction on $f$ is not enough. Recently, Muthuvel proved in \cite{MR3310341}, Theorem 3, that under the continuum hypothesis there exists an $\omega_1$-reasonable set mapping on $\omega_1$ with no infinite free set. He also proved that if the splitting number $\mathfrak{s}$ is above $\aleph_1$ (in which case, the continuum hypothesis fails) then for every $\omega_1$-reasonable set mapping on $\omega_1$ there is an infinite free set (\cite{MR3310341}, Theorem 1).

Muthuvel proved also that if $f:\mathbb{R}\rightarrow\mathcal{P}(\mathbb{R})$ satisfies $C(\omega,\omega)$ then there exists a dense free set for $f$ (in the usual topology), see \cite{MR3310341}, Corollary 1. This result is remarkable, but the proof required a very strong assumption. Unlike the assumptions $C(\kappa,\omega)$ and $C(\omega,\kappa)$ which are obligatory, the property $C(\omega,\omega)$ means that the range of the function is closed to be a collection of disjoint sets. The main objective of the current paper is to drop this assumption.

We need, however, to assume that the continuum hypothesis fails. We shall use the assumption that $\binom{\kappa}{\omega} \rightarrow \binom{\kappa}{\omega}^{1,1}_2$, when $\kappa=2^\omega$. It holds, e.g. if the reaping number $\mathfrak{r}$ is below the cofinality of the continuum. Under this assumption we will be able to show that every $\kappa$-reasonable function on the reals has an infinite dense free set.
In some sense, this theorem is the dual to the theorem of Muthuvel on $\omega_1$, since the reaping number $\mathfrak{r}$ is the dual of the splitting number $\mathfrak{s}$.
Let us recall the definition of these cardinal characteristics:

\begin{definition}
\label{rrrsss} The reaping and splitting numbers.
\begin{enumerate}
\item [$(\aleph)$] Suppose $B\in[\omega]^\omega$ and $S\subseteq\omega$. $S$ splits $B$ if $|S\cap B|=|(\omega\setminus S)\cap B|=\aleph_0$.
\item [$(\beth)$] $\{T_\alpha:\alpha<\kappa\}$ is an unreaped family if there is no $S\in[\omega]^\omega$ so that $S$ splits $T_\alpha$ for every $\alpha<\kappa$. Likewise, $\{S_\alpha:\alpha<\kappa\}$ is a splitting family iff for every $B\in[\omega]^\omega$ there exists an ordinal $\alpha<\kappa$ so that $S_\alpha$ splits $B$.
\item [$(\gimel)$] The reaping number $\mathfrak{r}$ is the minimal cardinality of an unreaped family, and the splitting number $\mathfrak{s}$ is the minimal cardinality of a splitting family.
\end{enumerate}
\end{definition}

An important tool in the proof of the main theorem below (as well as the proofs of Muthuvel) is the polarized partition realtion. It can be phrased in the language of colorings, or in the language of partitions (we shall use both):

\begin{definition}
\label{pppp} The strong polarized relation. \newline
We say that the strong polarized relation $\binom{\lambda}{\kappa} \rightarrow \binom{\lambda}{\kappa}^{1,1}_2$ holds iff for every coloring $c : \lambda \times \kappa \rightarrow 2$ there are $A \subseteq \lambda$ and $B \subseteq \kappa$ such that $|A|=\lambda, |B|=\kappa$ and $c \upharpoonright (A \times B)$ is constant.
\end{definition}

We try to use standard notation. We denote cardinals by $\theta,\kappa,\lambda,\mu$ and ordinals by $\alpha,\beta,\gamma,\delta$.
By $[A]^\theta$ we denote the collection of all subsets of $A$ of cardinality $\theta$. We use the Jerusalem forcing notation, so $p\leq q$ reads $q$ is stronger than $p$.

We mention the Erd\"os-Dushnik-Miller theorem which says that $\lambda\rightarrow(\lambda,\omega)^2$ for every infinite cardinal $\lambda$. A proof of this theorem appears in \cite{MR795592}.
For a general background about cardinal characteristics we refer to \cite{MR2768685}. For the above combinatorial theorems and related results we suggest \cite{MR795592} and \cite{MR3075383}. For basic account of forcing and Martin's axiom we advert to \cite{MR597342}.

I thank the referee for the careful reading, mathematical corrections and comments which improved the readability of this paper.

\newpage

\section{Dense free subsets of the reals}

We commence with the combinatorial theorem:

\begin{theorem}
\label{mt} Infinite free subsets. \newline
Suppose that $\kappa>\aleph_0, \binom{\kappa}{\omega} \rightarrow \binom{\kappa}{\omega}^{1,1}_2$ and $f$ is a $\kappa$-reasonable set mapping from $\kappa$ into $\mathcal{P}(\kappa)$. \newline
Then there exists an infinite free subset for $f$.
\end{theorem}

\par\noindent\emph{Proof}. \newline
We define a coloring $c: [\kappa]^2\rightarrow 2$ as follows. $c(\{\alpha,\beta\})=1$ iff $\alpha\notin f(\beta)\wedge \beta\notin f(\alpha)$. We employ the Erd\"os-Dushnik-Miller theorem to get either $H_0\in[\kappa]^\kappa$ such that $c\upharpoonright[H_0]^2=\{0\}$ or $H_1\in[\omega]^\omega$ such that $c\upharpoonright[H_1]^2=\{1\}$. If there exists such $H_1$ then we are done, since it would be a free set for $f$ by the definition of the coloring $c$, so assume towards contradiction that there is no $H_1$ as above.

We choose $A\in[H_0]^\omega,B\in[H_0]^\kappa$ so that $A\cap B=\emptyset$. We decompose the cartesian product $A\times B$ into two disjoint collections:

$$
A\times B=\{\langle a,b\rangle: a\in f(b)\}\bigcup \{\langle a,b\rangle: a\notin f(b)\}.
$$

By the assumption $\binom{\kappa}{\omega} \rightarrow \binom{\kappa}{\omega}^{1,1}_2$ we choose $A_0\in[A]^\omega,B_0\in[B]^\kappa$ such that either $A_0\times B_0\subseteq\{\langle a,b\rangle: a\in f(b)\}$ or $A_0\times B_0\subseteq\{\langle a,b\rangle: a\notin f(b)\}$ (here we use the language of partitions with respect to the polarized relation).

If $A_0\times B_0\subseteq\{\langle a,b\rangle: a\in f(b)\}$ then $A_0\subseteq \bigcap\{f(b):b\in B_0\}$, contradicting the assumption that $f$ is $C(\kappa,\omega)$. Similarly, if $A_0\times B_0\subseteq\{\langle a,b\rangle: a\notin f(b)\}$ then $b\in f(a)$ for every $a\in A_0,b\in B_0$ (since $a\notin f(b)$ and all the members are taken from $H_0$), so $B_0\subseteq \bigcap\{f(a):a\in A_0\}$, contradicting the assumption that $f$ is $C(\omega,\kappa)$.

\hfill \qedref{mt}

\begin{corollary}
\label{mmt} Assume $\mathfrak{r}<\cf(\kappa)\leq\kappa=2^\omega$. \newline
For every $\kappa$-reasonable function $f:\mathbb{R}\rightarrow\mathcal{P}(\mathbb{R})$ there exists an infinite free subset.
\end{corollary}

\par\noindent\emph{Proof}. \newline
By the assumption on the reaping number $\mathfrak{r}$, the relation $\binom{\kappa}{\omega} \rightarrow \binom{\kappa}{\omega}^{1,1}_2$ holds, as proved in \cite{MR3201820}, Claim 1.4. Hence the above theorem applies, and there exists an infinite free subset as desired.

\hfill \qedref{mmt}

\begin{remark}
\label{rremark}
\begin{enumerate}
\item [$(\alpha)$] If one wishes to assume stronger intersection properties for $f$, then the polarized partition relation can be weakened. Generally, the relation $\binom{\kappa}{\omega} \rightarrow \binom{\theta_0\ \theta_1}{\omega\ \omega}^{1,1}_2$ provides an infinite free subset for $f$ whenever $f$ is $C(\omega,\theta_0)$ and $C(\theta_1,\omega)$. The proof is just the same.
\item [$(\beta)$] Nevertheless, it is consistent that $\binom{\mu}{\omega} \nrightarrow \binom{\omega_1}{\omega}^{1,1}_2$ for every $\mu\in(\omega,2^\omega]$, by adding $\lambda$-many Cohen reals, see \cite{MR3000439}, Claim 2.4.
\item [$(\gamma)$] The strong polarized relation cannot be weakend here. Even if we assume, e.g. that $\binom{\omega_1}{\omega} \rightarrow \binom{\omega_1\ \alpha}{\omega\ \omega}^{1,1}_2$ for every $\alpha<\omega_1$, and $f$ is reasonable for $\omega_1$, we may fail to get an infinite free set. Indeed, this relation holds under the continuum hypothesis (see \cite{MR2444279}).
\item [$(\delta)$] It is tempting to try to generalize the above theorems upon replacing $\aleph_0$ by some larger cardinal. We can take a supercompat cardinal $\lambda$ in lieu of $\aleph_0$, and force for $\mu>\lambda$ that $\binom{\mu}{\lambda} \rightarrow \binom{\mu}{\lambda}^{1,1}_2$ as shown in \cite{MR2927607}. However, we need also the equivalent to the Erd\"os-Dushnik-Miller theorem, and this would give a monochromatic set only of size (arbitrarily large) less than $\lambda$.
\item [$(\varepsilon)$] The assumption $\mathfrak{r}<\cf(2^\omega)$ gives an infinite free set for reasonable functions $f:\theta\rightarrow\mathcal{P}(\theta)$ not only for the continuum but for every $\theta\geq\cf(\theta)>\mathfrak{r}$. Similarly, Theorem 1 of \cite{MR3310341} applies to every $\aleph_1\leq\theta<\mathfrak{s}$ so that $\cf(\theta)>\aleph_0$. It follows that one can prove the consistency of free sets for each $\theta\in[\aleph_1, 2^{\aleph_0}]$ simultaneously, see \cite{MR3201820}.
\end{enumerate}
\end{remark}

We turn now to the topological density of the free set. We focus on a set mapping defined on the reals, and we are looking for a free set which is also dense. The existence problem of such sets appears in \cite{MR1086127}.
Our proof is just as in \cite{MR3310341}, but we can replace the property $C(\omega,\omega_1)$ of the function $F$ (defined below) by the weaker demand $C(\omega,\kappa)$.

\begin{theorem}
\label{mt1} Dense free set of the reals. \newline
Let $\kappa$ be $2^\omega$, and assume $f:\mathbb{R}\rightarrow\mathcal{P}(\mathbb{R})$ is $\kappa$-reasonable (i.e., satisfies $C(\kappa,\omega)$ and $C(\omega,\kappa)$). \newline
If $\binom{\kappa}{\omega} \rightarrow \binom{\kappa}{\omega}^{1,1}_2$ then there exists an everywhere dense free set for $f$.
\end{theorem}

\par\noindent\emph{Proof}. \newline
Firstly, we define a derived function $F:\mathbb{R}\rightarrow\mathcal{P}(\mathbb{R})$ as follows:

$$
F(x)=\{y:y\in f(x) \vee x\in f(y)\}.
$$

We claim that $F$ is $C(\omega,\kappa)$. For proving this fact, assume towards contradiction that there exists a subset $A\subseteq\mathbb{R}, |A|=\aleph_0$ such that $D=\bigcap\{F(x):x\in A\}$ is of size $\kappa$. By removing a countable subset from $D$ we may assume without loss of generality that $A\cap D=\emptyset$. We write:

$$
A\times D=\{\langle a,b\rangle:a\in f(b)\}\cup \{\langle a,b\rangle:a\notin f(b)\}.
$$

Since $\binom{\kappa}{\omega} \rightarrow \binom{\kappa}{\omega}^{1,1}_2$ we can choose $H_0\in[A]^\omega, H_1\in[D]^\kappa$ so that $H_0\times H_1\subseteq \{\langle a,b\rangle:a\in f(b)\}$ or $H_0\times H_1\subseteq \{\langle a,b\rangle:a\notin f(b)\}$. If $H_0\times H_1\subseteq \{\langle a,b\rangle:a\in f(b)\}$ then $H_0\subseteq\bigcap\{f(b): b\in H_1\}$, contradicting the assumption that $f$ is $C(\kappa,\omega)$. Likewise, if $H_0\times H_1\subseteq \{\langle a,b\rangle:a\notin f(b)\}$ then $b\in f(a)$ for every $\langle a,b\rangle\in H_0\times H_1$ and hence $H_1\subseteq \bigcap\{f(a): a\in H_0\}$, contradicting the assumption that $f$ is $C(\omega,\kappa)$.

Having proved that $F$ is $C(\omega,\kappa)$ we can build a dense free set $S=\{x_m:m\in\omega\}$ by induction on $m\in\omega$. Let $\{I_n:n\in\omega\}$ enumerate all the finite open intervals of the reals with endpoints from $\mathbb{Q}$. We choose a subset $E_0\subseteq\mathbb{R}$ of size $\kappa$ such that $|E_0\cap I_n|=\kappa$ for every $n\in\omega$. Along the induction, we keep the fact that $E_m\cap I_n$ is of size $\kappa$, for every $m,n\in\omega$.
We describe the choice of the first member $x_0$ of the set $S$. Define the set $C^0_n$ for every $n\in\omega$ as follows:

$$
C^0_n=\{x\in E_0\cap I_0:|(E_0\setminus F(x))\cap I_n|<\kappa\}.
$$

We wish to prove that $C^0_n$ is finite (for every $n\in\omega$), so fix a natural number $n$ and assume to the contrary that $C^0_n$ is infinite. We choose an infinite countable subset $C\subseteq C^0_n$, and we get the bound $|\bigcup\{(E_0\setminus F(\ell))\cap I_n:\ell\in C\}|<\kappa$, being a coubtable union of sets of size less than $\kappa$ (notice that $\cf(\kappa)>\aleph_0$). However, $F$ is $C(\omega,\kappa)$, so the cardinality of $\bigcap\{F(\ell):\ell\in C\}$ is less than $\kappa$. Since the cardinality of $E_0\cap I_n$ equals $\kappa$, we have:

$$
|(E_0\cap I_n) \setminus \bigcap\{F(\ell):\ell\in C\}|=\kappa.
$$

This fact leads to a contradiction. Indeed, $(E_0\cap I_n) \setminus \bigcap\{F(\ell):\ell\in C\} = (E_0\setminus \bigcap\{F(\ell):\ell\in C\})\cap I_n = \bigcup \{E_0\setminus F(\ell):\ell\in C\}\cap I_n = \bigcup\{(E_0\setminus F(\ell))\cap I_n:\ell\in C\}$, and the latter is of size less than $\kappa$ as we have seen before.

Concluding that $C^0_n$ is finite we infer that the size of $\bigcup\{C^0_n:n\in\omega\}$ is countable, so we can choose $x_0\in (E_0\cap I_0)$ such that $x_0\notin \bigcup\{C^0_n:n\in\omega\}$. This is the first step of the induction.

Suppose $x_m$ is at hand, and the sets $C^m_n=\{x\in E_m\cap I_m:|(E_m\setminus F(x))\cap I_n|<\kappa$ were defined similarly to $C^0_n$ and satisfy $|C^m_n|<\aleph_0$ for every $n\in\omega$.
We need to choose $x_{m+1}$. Let $E_{m+1}$ be $E_m\setminus F(x_m)$. Notice that $|E_{m+1}|=\kappa$, moreover $|E_{m+1}\cap I_n|=\kappa$ for every $n\in\omega$, by the properties of $C^m_n$. Hence we can pick up $x_{m+1}\in E_{m+1}\cap I_{m+1}$ so that $x_{m+1}\notin\{x_j:j<m+1\}$.

Finally, let $S=\{x_m:m\in\omega\}$. Clearly, $S$ is a dense set in the usual topology of the reals, as $S$ intersects every open interval. Likewise, $S$ is a free set for $f$. Indeed, if $i<j<\omega$ then $x_j\notin F(x_i)$ by the construction, hence $x_j\notin f(x_i)$ and $x_i\notin f(x_j)$ by the definition of $F$, so we are done.

\hfill \qedref{mt1}

It may help to notice that there are $\kappa$-many countable dense free subsets for every $f$ with the assumed properties, since at each stage one can choose $x_m$ from a set of size $\kappa$.

Our last theorem generalizes Theorem 2 of \cite{MR3310341}, where he proves that under Martin's axiom one can get uncountable free sets for reasonable functions on $\omega_1$. We shall see that, under MA + $\neg{\rm CH}$, for a regular uncountable cardinal $\kappa<2^{\aleph_0}$, one can get free sets of size $\kappa$ for every $f:\kappa\rightarrow\mathcal{P}(\kappa)$ which satisfies $C(\omega,\omega_1)$ and $C(\omega_1,\omega)$. Since we are confined to $ccc$ forcing notions, we cannot weaken the assumption on $f$ into $C(\omega,\kappa)$ and $C(\kappa,\omega)$. However, we introduce the forcing notion which gives a free set for a specific $f$ under the weak assumption. The problem is to iterate such forcings in order to cover all possible functions.

\begin{definition}
\label{fffnotions} The forcing notions $\mathbb{P}_f$ and $\mathbb{P}_f^{\rm fin}$. \newline
Assume $\omega_1\leq\kappa=\cf(\kappa)\leq 2^\omega$ and $f:\kappa\rightarrow\mathcal{P}(\kappa)$.
\begin{enumerate}
\item [$(\aleph)$] $p\in\mathbb{P}_f$ iff $p\subseteq\kappa,|p|<\kappa$ and $p$ is free for $f$.
\item [$(\beth)$] $p\in\mathbb{P}_f^{\rm fin}$ iff $p\subseteq\kappa,|p|<\aleph_0$ and $p$ is free for $f$.
\item [$(\gimel)$] The order (in both notions) is $p\leq q$ iff $p\subseteq q$.
\end{enumerate}
\end{definition}

The generic object of both forcing notions gives a free set of size $\kappa$ for the function $f$. The following claim shows that the above forcing notions preserve cardinals and cofinalities. The advantage of $\mathbb{P}_f^{\rm fin}$ is the $ccc$, but the price is the assumption that $f$ is $\omega_1$-reasonable rather than $\kappa$-reasonable.

\begin{claim}
\label{nnnotcollapsing} Assume $\omega_1\leq\kappa=\cf(\kappa)< 2^\omega, f:\kappa\rightarrow\mathcal{P}(\kappa)$ and $\mathfrak{p}=2^\omega$.
\begin{enumerate}
\item [$(a)$] If $f$ is $\kappa$-reasonable then $\mathbb{P}_f$ is $\kappa$-complete and $\kappa^+$-cc.
\item [$(b)$] If $f$ is $\omega_1$-reasonable then $\mathbb{P}_f^{\rm fin}$ is $ccc$.
\end{enumerate}
\end{claim}

\par\noindent\emph{Proof}. \newline
The completeness assertion in $(a)$ follows by the fact that $\kappa$ is a regular cardinal. For the chain condition the proof of $(a)$ and $(b)$ is essentialy the same, so we focus on the $ccc$ of part $(b)$, and we indicate that for $(a)$ one needs to replace $\omega_1$ by $\kappa^+$.

Assume towards contradiction that $\mathcal{A}=\{p_\alpha:\alpha<\omega_1\} \subseteq \mathbb{P}_f^{\rm fin}$ is an antichain.
It means that for every $\alpha<\beta<\omega_1$ there are $x\in p_\alpha, y\in p_\beta$ so that $x\in f(y)\vee y\in f(x)$.
By the Delta-system lemma we can shrink $\mathcal{A}$ into a set of size $\omega_1$ for which every pair of members has the same intersection $r$. Since the members of $r$ are free with all the members of the conditions in $\mathcal{A}$, we assume without loss of generality that $r=\emptyset$. Moreover, we assume without loss of generality that $|p_\alpha|=n$ for every $p_\alpha\in\mathcal{A}$, and we write $p_\alpha=\{x_i^\alpha:i<n\}$. We choose a nonprincipal ultrafilter $U$ on $\omega$.

Since $\mathfrak{p}=2^\omega$, $\mathfrak{s}=2^\omega$ as well and hence $\binom{\kappa}{\omega} \rightarrow \binom{\kappa}{\omega}^{1,1}_2$ for every $\kappa\in(\omega,2^\omega)$ with uncountable cofinality, including $\omega_1$.
In order to use the assumed polarized relation, we decompose $\mathcal{A}$ into two disjoint sets $\{p_\ell:\ell\in\omega\}\cup \{p_\alpha:\omega\leq\alpha<\omega_1\}$. Fix an ordinal $\omega\leq\alpha<\omega_1$. By the assumption towards contradiction, for every $\ell\in\omega$ there is a member of $p_\ell$ which is not free with a member of $p_\alpha$, and hence:

$$
\bigcup\limits_{i,j<n}\{\ell\in\omega:x^\ell_i\in f(x^\alpha_j)\vee x^\alpha_j\in f(x^\ell_i)\}=\omega\in U.
$$

The above is a finite union, so for some $i,j<n$ we have $A^{ij}_\alpha=\{\ell\in\omega:x^\ell_i\in f(x^\alpha_j)\vee x^\alpha_j\in f(x^\ell_i)\}\in U$. Notice that $i,j$ depend on $\alpha$, and $\{A^{ij}_\alpha: \omega\leq\alpha<\omega_1\}\subseteq U$. Since $\mathfrak{p}=2^\omega$ we can find $A\in[\omega]^\omega,B\in[\omega_1]^{\omega_1}$ so that $A\subseteq A^{ij}_\alpha$ for every $\alpha\in B$. Without loss of generality, the ordinals $i,j$ are the same for every $\alpha\in B$, by shrinking $B$ again if needed.

The cartesian product $A\times B$ is expressible now as $\{(\ell,\alpha): x^\ell_i\in f(x^\alpha_j)\}\cup\{(\ell,\alpha): x^\ell_i\notin f(x^\alpha_j)\}$. Since $\binom{\omega_1}{\omega} \rightarrow \binom{\omega_1}{\omega}^{1,1}_2$ we can find $A_0\in[A]^\omega,B_0\in[B]^{\omega_1}$ such that $A_0\times B_0 \subseteq \{(\ell,\alpha): x^\ell_i\in f(x^\alpha_j)\}$ or $A_0\times B_0 \subseteq \{(\ell,\alpha): x^\ell_i\notin f(x^\alpha_j)\}$. But this contradicts either that $f$ is $C(\omega,\omega_1)$ or that $f$ is $C(\omega,\omega_1)$, so we are done.

\hfill \qedref{nnnotcollapsing}

Having established the chain condition, we can iterate in order to create a free set for every $f$. Notice that $\mathbb{P}_f$ preserves cardinals but the iteration is more involved, so we iterate $\mathbb{P}_f^{fin}$.
We employ Martin's axiom, although one needs only the combinatorial assumptions of the above claim.

\begin{theorem}
\label{mt2} Martin's axiom and $\kappa$-free sets. \newline
Assume Martin's axiom and $2^{\aleph_0}>\aleph_1$. \newline
Then for every $\aleph_0<\kappa=\cf(\kappa)<2^{\aleph_0}$, and every $f:\kappa\rightarrow\mathcal{P}(\kappa)$ which is $\omega_1$-reasonable, there exists a free set of size $\kappa$.
\end{theorem}

\par\noindent\emph{Proof}. \newline
Under Martin's axiom we have $\mathfrak{p}=2^\omega$.
Given any such function $f$, we know that $\mathbb{P}_f^{\rm fin}$ is $ccc$. By Martin's axiom we can choose a generic set $G\subseteq \mathbb{P}_f^{\rm fin}$. We define $T=\bigcup G$. The natural dense subsets show that $T$ is unbounded in $\kappa$, hence $|T|=\kappa$ by the regularity of $\kappa$. Since $G$ is directed we infer that $T$ is a free set.

\hfill \qedref{mt2}

We conclude with several open problems raised by the above theorems. First, one may wonder if $\binom{\kappa}{\omega} \rightarrow \binom{\kappa}{\omega}^{1,1}_2$ is the correct assumption for large free sets to exist:

\begin{question}
\label{q0} Negative partition relation and free sets. \newline
Assume $\binom{\kappa}{\omega} \nrightarrow \binom{\kappa}{\omega}^{1,1}_2$. Can we find a $\kappa$-reasonable function on $\kappa$ with no free set of size $\kappa$ (or even without an infinite free set)?
\end{question}

\begin{question}
\label{q1} Cohen reals and free sets. \newline
Suppose $f$ is a reasonable function on the reals, in the universe forced by adding $\lambda$-many Cohen reals. Is it true that $f$ has an infinite free set?
Under the same assumption, is it true that $f$ has a dense free set?
\end{question}

Second, we can ask about stronger assertions to be forced:

\begin{question}
\label{q2} Free sets of size $2^\omega$. \newline
suppose $\kappa=2^\omega$. Is it consistent that every $\kappa$-reasonable $f$ from $\kappa$ into $\mathcal{P}(\kappa)$ has a free set of size $\kappa$?
\end{question}

\begin{question}
\label{q3} Weaker intersection properties. \newline
Can we replace the assumption that $f$ is $\omega_1$-reasonable by the weaker assumption that $f$ is $\kappa$-reasonable in Theorem \ref{mt2}?
\end{question}

For the last questions it seems that Martin's axiom is not sufficent.

\newpage

\bibliographystyle{amsplain}
\bibliography{arlist}

\end{document}